\documentclass[a4paper,aps,prl,onecolumn,notitlepage]{revtex4-1}

\usepackage{wrapfig}
\usepackage[small,bf,justification=raggedright]{caption}
\usepackage{graphicx} %[dvips]
\usepackage[usenames,dvipsnames]{xcolor}
\usepackage{tikz}
\usepackage{amsmath}
\usepackage{amssymb}
\usepackage{color}
\usepackage{natbib}
\begin{document}
\title{{\em Germain Curvature}: The Case for Naming the Mean Curvature of a Surface after Sophie Germain}
\author{Douglas P. Holmes$^1$}
\affiliation{\footnotesize$^1$ Mechanical Engineering, Boston University, Boston, MA, 02215, USA} 

\date{\today}

\maketitle
How do we characterize the shape of a surface? It is now well understood that the shape of a surface is determined by measuring how curved it is at each point. From these measurements, one can identify the directions of largest and smallest curvature, {\em i.e.} the principal curvatures, and construct two surface {\em invariants} by taking the average and the product of the principal curvatures. The product of the principal curvatures describes the {\em intrinsic} curvature of a surface, and has profound importance in differential geometry -- evidenced by the Gauss's {\em Theorema Egregium} and the Gauss--Bonnet theorem. This curvature is commonly referred to as the {\em Gaussian Curvature} after Carl Friedrich Gauss, following his significant contributions to the emerging field of differential geometry in his 1828 work~\cite{Gauss1828}. The average, or {\em mean curvature}, is an {\em extrinsic} measure of the shape of a surface -- that is, the shape must be {\em embedded} in a higher dimensional space to be measured. As the analogy goes, an ant living on a sheet of paper would not be able to know if that sheet was flat or rolled into a cylindrical shape without viewing the surface from afar. Since Gauss's {\em Theorema Egregium} states that the intrinsic curvature cannot change without stretching the surface, the extrinsic curvature is an important measure when describing the deformations of a surface in the absence of stretching -- like those often observed in slender structures. Beginning in 1811, and culminating in efforts in 1821 and 1826, the mathematician Sophie Germain identified the mean curvature as the appropriate measure for describing the shape of vibrating plates. Her hypothesis leads directly to the equations describing the behavior of thin, elastic plates. In letters to Gauss, she described her notion of a ``sphere of mean curvature'' that can be identified at each point on the surface. This contribution stimulated a period of rapid development in elasticity and geometry, and yet Germain has not yet received due credit for introducing the concept of mean curvature. It is clear from the primary source evidence described below that this measure of {\em mean curvature} should bear the name of Sophie Germain.

\subsection{Prix Extraordinaire}

\begin{wrapfigure}{r}{0.25\columnwidth}
\begin{center}
	\vspace{-8mm}
		\resizebox{0.25\columnwidth}{!} {\includegraphics{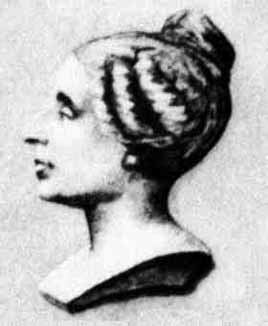}}
\end{center}\vspace{-7mm}
	\caption[]{Sophie Germain
	\label{portrait}
}
\vspace{-6mm}
\end{wrapfigure}
In 1808, Ernst Chladni showed a remarkable experiment to Pierre Simon Laplace. Chladni took a thin plate, covered its surface with powder, and vibrated it by running a bow along its edge. The grains of powder rose and fell following the plate's vibrations, coming to rest in its nodes. Each and every musical note, or frequency, played on the plate caused the grains to come to rest in different patterns. So too did changing its shape. Laplace, intrigued by these striking experiments and seeing an opportunity to elevate the stature of his prot\'{e}g\'{e}, Sim\'{e}on Denis Poisson, introduced Chladni and his experiment to Napoleon Bonaparte. Napoleon urged the First Class of the Institute of France to set the problem as the subject of a {\em prix extraordinaire}. The competition was announced in 1809, with 3000 francs to be awarded to the scholar who could provide a mathematical analysis to explain the vibration of thin, elastic plates. The prize was eventually won by Sophie Germain (1776--1831), a self-taught mathematician from Paris, France.

Entries to the competition were due in 1811. As members of the Institute, Laplace, Joseph-Louis Lagrange, and Adrien-Marie Legendre were to serve as judges. Mathematicians who will later play essential roles in the development of elasticity theory -- Poisson, Joseph Fourier, Claude-Louis Navier, and Augustin-Louis Cauchy -- did not enter the competition. Germain's contribution was the only submission. Germain began with an ``ingenious hypothesis'' that she developed by drawing an analogy to Leonhard Euler's model for a vibrating beam~\cite{Bucciarelli1980}. She anticipated that the elastic force in the plate would be proportional to the average of its two principal curvatures, its {\em mean curvature}. As we will see, this hypothesis both challenged the prevailing mentality of the time and intuited aspects of the geometry of surfaces that were not to be fully developed for decades. Although the equation Germain derived from this hypothesis was incorrect, the hypothesis itself leads directly to correct equation, 
\begin{equation}
\label{plateEQ}
\ddot{w} = k^2\nabla^4w,
\end{equation}
following the correct application of the calculus of variations~\footnote{Here, $w(x,y)$ is the out-of-plane deflection of the plate, $k^2$ is a ratio between the plate's bending rigidity and its mass (typically given by its areal density times its thickness), and $\nabla^4(\cdot)=\nabla^2\nabla^2(\cdot)$ is the biharmonic operator.}. Lagrange himself, while judging the contest, presented this as the correct equation {\em if} Germain's hypothesis was correct. At this point, there was little reason to believe in either Germain's hypothesis or the equation Lagrange derived from it, and so no prize was awarded.

The competition was extended until 1813. Legendre, who was in correspondence with Germain at the time~\cite{Bucciarelli1980}, shared Lagrange's equation with Germain along with advice on how to derive it: ``\ldots follow Lagrange's method in the new edition ({\em Mecanique Analytique}), page 148, adjusting the appropriate term representing the force due to elasticity'' (Legendre to Germain, Dec. 1811)~\cite{Stupuy1896, Bucciarelli1980}. Germain was unable to take this advice. This gets to the heart of the matter: ``\ldots all the evidence argues that Sophie Germain had a mathematical brilliance that never reached fruition due to a lack of rigorous training available only to men.'' Germain was familiar with Lagrange's {\em Mecanique Analytique}, but only through notes she was able to acquire from students who had studied under Lagrange at Ecole Polytechnique. She worked on her elasticity problem largely in isolation, save for some intermittent correspondence with Legendre. She struggled to apply Lagrange's variational calculus to the double integrals that appear in the formulation of a plate, a fact that can certainly be attributed to the academic barriers and isolation Germain faced due to her gender. 

%And yet, as her biographers Louis Bucciarelli and Nancy Dworsky make clear, ``\ldots even if she had become superbly proficient in variational techniques and thoroughly at home with physical principles, there is no likelihood that this would have qualified her for entrance into the scientific community'' ~\cite{Bucciarelli1980}.

Her second attempt at the prize went beyond an attempt to recover equation~\ref{plateEQ}. She identified some of the relevant boundary conditions, and used the equation Lagrange had derived from her hypothesis to accurately predict some node shapes and frequencies of square and rectangular plates. However, once again her effort yielded an incorrect derivation of the governing equation, and her justification of her hypothesis did not satisfy the judges~\footnote{Lagrange, who may have been sympathetic to her treatment of the problem, passed away in April 1813. Entries were due in October 1813.}. As we will see, accepting her hypothesis would have required a substantial shift in thinking for the time. Nevertheless, the accurate node shape predictions made her hypothesis even more intriguing, and while no prize was awarded, Germain received an honorable mention, and the prize was extended a second time.

Her third attempt at the prize extended her hypothesis to shallow shells, contending that the elastic force was proportional to the ``\ldots difference in the shape of the deformed and undeformed surfaces'', {\em i.e.} the change in the mean curvature of the surface. She complemented these efforts with experimental work on vibrations of shallow, cylindrical shells. Despite not finding a justification of her hypothesis that was suitable to the committee, they were sufficiently impressed with her solutions to particular integrals, along with her experimental efforts, and she was awarded the prize in January 1816.
%Her difficulties in employing the variational methods resulted in a theoretical portion of her report that was ``\ldots interesting in intent, but fundamentally deficient and wrong in execution''~\cite{Bucciarelli1980}. 

To understand the significance of Germain's hypothesis, it is essential to understand two things. First, it was in direct opposition to the {\em molecular mentality}, a hypothesis that was firmly held by most scholars during this time -- including Laplace and Poisson, who were judging her work.  Second, the differential geometry of surfaces was not yet established. Taken together, Germain's work makes an unassailable case for referring to the {\em mean curvature} as the {\em Germain curvature}.

\begin{figure}[t]
\includegraphics[width=\columnwidth]{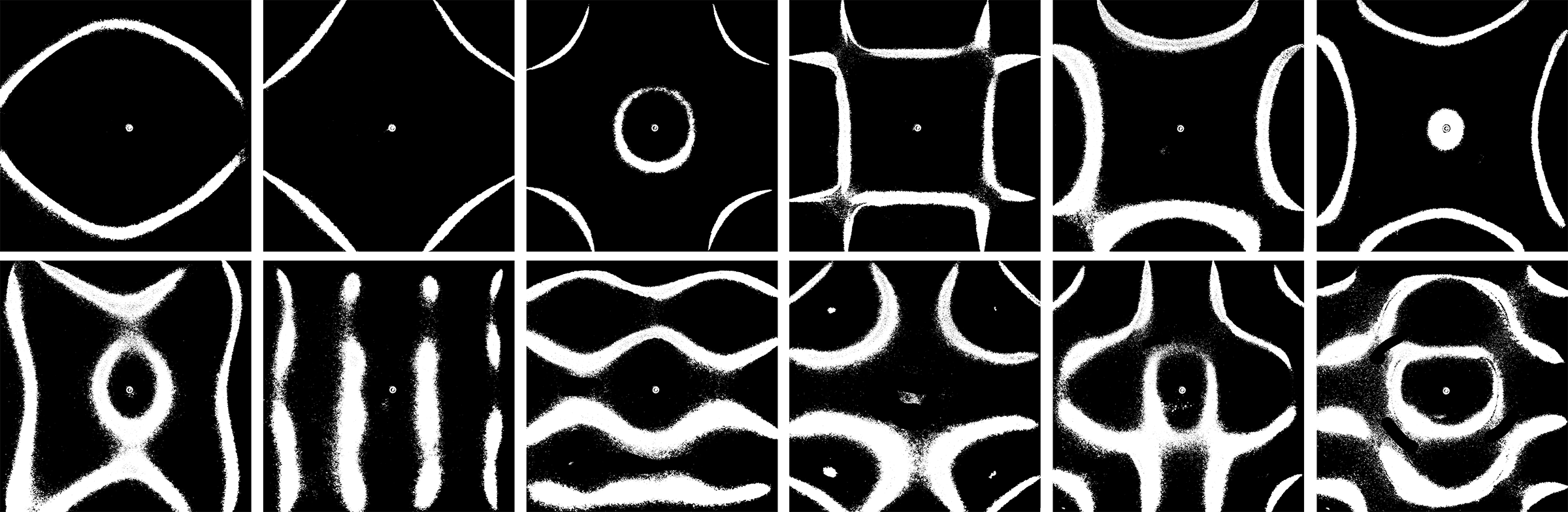}
\caption[]{Experimental images of Chladni patterns formed by vibrating a square plate (30.5mm x 30.5mm x 0.37mm) at different frequencies. Frequencies range from 100 Hz (Top Left) to 650 Hz (Bottom Right), incrementing by 50 Hz in each frame. Salt is added to the surface to illustrate the nodes. \label{Chladni}}
\end{figure}

\subsection{Molecular Mentality}
The molecular mentality presumed that there must be ``sensible forces at insensible distances,'' and thus it was believed that Germain's hypothesis, which was based on the deformations of an elastic surface, was simply an intermediate step on the way to a theory based on ``molecular actions''. Poisson was unconvinced by the variational methods introduced by Lagrange in his {\em Mecanique Analytique}, and set out to derive equation~\ref{plateEQ} from molecular principles.  As a result, Poisson derived an incorrect, nonlinear equation to describe the bending of a plate from a model that does not include bending~\cite{Bucciarelli1980}. Germain had full confidence in her hypothesis, believing that ``whatever may be going on inside the surface, elastic behavior is measured by deformation.''~\cite{Bucciarelli1980}, and it was argued that her work ``\ldots despite its flaws, became one (of several) rallying points for criticism of that molecular mentality''~\cite{Bucciarelli1980}. Germain's efforts caught the attention of Fourier and Navier, and later Cauchy, who each had their own approach to the problem. In 1818, Fourier demonstrated how a solution to the plate equation could be obtained using the mathematical method that now bears his name -- that is, he represented the plate's motion as the sum over trigonometric functions. His efforts focused on the solution of a specific problem: the response of an {\em infinite} plate to an initial disturbance~\footnote{This problem enabled him to avoid one of the problems in the theory at that time: a general, accepted set of boundary conditions had not yet been established. By working with an infinite plate, Fourier was able to avoid the question of what was happening at the boundaries altogether.}. Navier took a more practical approach to understand plate bending. He considered the static loading of a plate supported at its edges. In the opening of his memoir, Navier acknowledged the importance of Germain's work (emphasis mine):

\begin{quote}
``The interesting experiments of M. Chladni on the vibrations of plates have stimulated the application of the calculus to the laws of movement what are manifest in these experiments: this was the subject of a prize proposed by the First Class of the Institute and won by Mademoiselle Germain. The research that was awarded {\bf the prize was founded on an ingenious hypothesis}, namely, that flexure gave birth, at each point of an elastic plate, to a force proportional to the sum of the inverse values of the two radii of principal curvature. Mademoiselle Germain gave the differential equation of equilibrium and movement of an elastic plane and some integrals of these equations, analogous to those Euler has given for the elastic lamina.'' -- Navier~\cite{Navier1820, Bucciarelli1980}
\end{quote}

Regardless of the respect that both Navier and Fourier had for Germain~\footnote{Fourier and Germain's friendship began in 1816 and lasted over a decade. See Ref.[2], p.87.}, they still maintained aspects of the molecular mentality in their efforts to describe the behavior of elastic plates. These efforts stimulated Cauchy's interest in the problem~\footnote{Cauchy was tasked with reviewing Navier's works -- three memoirs over two years -- submitted for consideration to the Academy of the Institute. Navier felt Cauchy procrastinated on reviewing his work to further Cauchy's own efforts in establishing a theory of elasticity. Later, Poisson presented his own memoir on the topic, omitting any mention of Navier. Navier took strong objection to both, claiming priority on the development of the theory of elasticity. Ref. [2], Ch.9.}, who approached the problem from a continuum perspective. In fact, in his critique of Navier's work, Cauchy introduced a definition of stress, a derivation of the equilibrium equations, a definition of strain, a derivation of the kinematic equations, and suggested several constitutive models~\cite{Bucciarelli1980, Cauchy1823, Truesdell1992}. At this point, it is clear that mathematicians had moved beyond the specific problem of a vibrating plate to the more general behavior of arbitrarily shaped elastic solids. Germain had been awarded the {\em prix extraordinaire}, but she remained outside the scientific community, and largely isolated beyond her friendship with Fourier. As her biographers note,
\begin{quote}
``The elasticity problem had now moved into a community to which she did not and could not belong, even though she did not know it, and even though that community would not or could not really tell her so. This difficulty existed quite independently of her mathematical inadequacies, though these inadequacies came from the same source, namely, her being a woman and thus excluded from formal education. \ldots When elasticity became a topic of interest among the professionals, a woman had no place in their midst.'' -- L. Bucciarelli and N. Dworsky~\cite{Bucciarelli1980}.
\end{quote}

Accepting the notion that a plate's behavior can largely be determined by its deformation required a break from the molecular mentality. Measuring that deformation required an understanding of the geometry of curved surfaces, a branch of mathematics that was developing in parallel at the time. 

Consider the challenge of measuring how a plate is bent~\cite{Bucciarelli1980}:
\begin{quote}
``To establish an expression for the curvature of a surface was not a simple problem. In the case of the beam, which Sophie Germain cited in support of her argument, there is little difficulty; the deformed beam is represented by one curve in space, with one radius of curvature. The elastic deformed surface, however, presents a multitude of possible curves through any point on the surface.'' 
\end{quote}
\subsection{Geometry of Surfaces}
Indeed, the reluctance to accept Germain's notion that the principal curvatures give the measure of curvature of a surface at a point was one of Poisson's primary objections to her hypothesis. Poisson was referring to Euler's ``\ldots discussion of the infinite number of different curves obtained from the intersections of different planes passing through a given point of the surface''~\cite{Bucciarelli1980}. Germain had studied Euler's works~\cite{Bucciarelli1980}, and was thus familiar with the method's Euler developed for calculating the curvatures through a point on a surface~\cite{Euler1767,Bardini2016}. Germain had intuited the significance of the principal curvatures in her first entry to the {\em prix extraordinaire}, in which she justified her choice of the mean curvature by noting that ``\ldots $(1/r)\cdot(1/r')$, for example, may always be neglected with respect to $((1/r)+(1/r')$''~\cite{Bucciarelli1980}\footnote{Here, $r$ and $r'$ are the principal radii of curvature.}.

\begin{wrapfigure}{r}{0.28\columnwidth}
\begin{center}
	\vspace{-5mm}
		\resizebox{0.28\columnwidth}{!} {\includegraphics{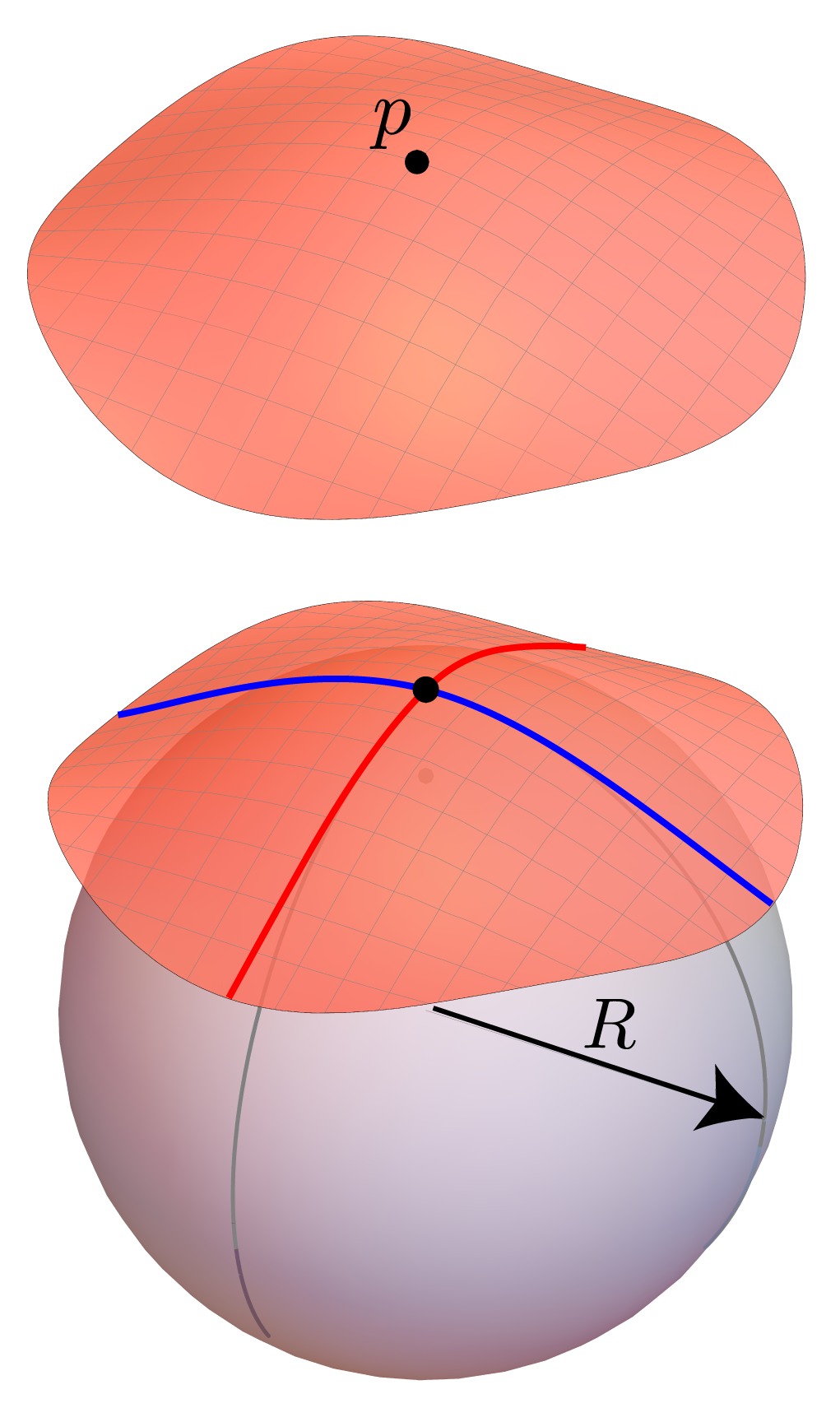}}
\end{center}\vspace{-7mm}
	\caption[]{A {\em sphere of mean curvature} at a point $p$ on a surface. 
	\label{portrait}
}
\vspace{-13mm}
\end{wrapfigure}
The elephant in the room is, of course, Carl Friedrich Gauss, whose name is now ascribed to the first measure of curvature Germain mentioned above, {\em i.e.} $(1/r)\cdot(1/r')$ is the {\em intrinsic curvature} or  {\em Gaussian curvature}. Prior to her work on elasticity, Germain was drawn to number theory after reading Gauss's {\em Disquisitiones Arithmeticae}. She adopted the pen name of Antoine-August LeBlanc, a student attending Ecole Polytechnique, and began corresponding with Gauss in 1804~\cite{Del2012}. The content of these letters reveal a mathematician who was eager to correspond with someone about his ``higher arithmetic''. While Germain always responded immediately to Gauss's letters, Gauss responded promptly on only one occasion: upon learning he had been corresponding with a woman~\cite{Bucciarelli1980, Mackinnon1990}. Here again, Napoleon plays a pivotal role in this story. In 1806, as Napoleon's armies engaged with Prussian forces, Germain -- recalling the fate of Archimedes~\cite{Dijksterhuis2014}\footnote{The account of Archimedes' death that Germain was aware of suggested that Archimedes was deeply engaged with a mathematical problem when the city was captured. It is believed that a Roman soldier commanded him to come and meet Marcellus, but he declined, preferring to finish working on the problem. The soldier killed Archimedes with his sword. See Ref.~\cite{Dijksterhuis2014}.} -- requested a favor from a family friend who was commanding the French artillery in Prussia: discover Gauss's whereabouts and ensure he remained safe~\footnote{A battalion commander traveled two hundred miles west of where he was stationed to check on Gauss.}. Gauss was confused, as he did not know a woman named Germain in Paris. Germain acknowledged ``\ldots that fearing the ridicule attached to a female scientist, I have previously taken the name of M. LeBlanc in communicating to you those notes that, no doubt, do not deserve the indulgence with which you have responded''~\cite{Bucciarelli1980}. Gauss responded immediately, expressing sincere appreciation for Germain's actions to ensure his safety, and respect and admiration for her mathematical curiosity and ability~\cite{Mackinnon1990}, 
\begin{quote}
``\ldots when a woman, because of her sex, our customs and prejudices, encounters infinitely more obstacles than men in familiarizing herself with their knotty problems, yet overcomes these fetters and penetrates that which is most hidden, she doubtless has the most noble courage, extraordinary talent, and superior genius.'' -- Gauss to Germain, April 1807~\cite{Bucciarelli1980}.
\end{quote}

Gauss and Germain remained in correspondence for the next couple of years, and faded as Germain's focus shifted from number theory to elasticity. It remains an important part of this story though, as historians believe that Gauss was likely made aware of Germain's work on Chladni's patterns and her measure of curvature~\cite{Bucciarelli1980}. As Germain still wrestled with how to justify her hypothesis, she produced memoirs in 1821 and 1826 that went into greater detail about measuring the curvature of a surface at a point. In these works, she introduced the notion of a reference sphere, whereby the radius of a sphere defines the mean curvature of the surface at a point in the same way that the radius of a reference circle defines the curvature at a point along a curve. Germain learned of Gauss's treatise on geometry~\cite{Gauss1828} in early 1829, when Gauss advised his student, M. Bader, to meet Germain during his visit to Paris.
\begin{quote}
``In chatting with Monsieur Bader about the current subject of my study, I provided him with the occasion to speak to me, and subsequently, to show me, the learned memoir in which you compare the curvature of surfaces to that of a sphere. \ldots I cannot tell you, Monsieur, how astonished, and at the same time, how satisfied I was in learning that a renowned mathematician, almost simultaneously, had the idea of an analogy that seems to me so rational that I neither understood how no one had thought of it sooner, nor how no one has wished to give any attention to date to what I have already published in this regard.'' -- Germain to Gauss, March 1829~\cite{Bucciarelli1980}.
\end{quote}
In the remainder of the letter, Germain illustrates her firm understanding of this measure of curvature. In rebutting Poisson's objection to her hypothesis, she remarked to Gauss that when
\begin{quote}
``\ldots taking the sum of the inverses of the radii of curvature of two of the curves of normal intersection, it would suffice to choose those that are contained in two mutually perpendicular planes in order to obtain a constant quantity equal to the sum of the inverses of the two radii of principal curvature. Furthermore, I remarked that if the two planes of intersection chosen make an angle of 45 degrees with those that contain the lines of principal curvature, the radii of curvature of the curves contained in the chosen planes are equal to each other and to the radius of a sphere that cuts the surface in such a way that two of the quadrants will be blow and the two others above the surface. This sphere is what I have named {\em the sphere of mean curvature}.'' -- Germain to Gauss, March 1829~\cite{Bucciarelli1980}\footnote{Her geometric description is reminiscent of the graphical tools developed by Karl Culmann and Christian Otto Mohr, now bearing the name of ``Mohr's Circle'', which aid in the transformation of any symmetric 2x2 tensor matrix, such as the curvatures of a surface.} .
\end{quote}

It is this notion of a {\em sphere of mean curvature} that Germain clearly introduced that should bear her name. 

\subsection{Germain Curvature}

It seems that every few decades the story of Sophie Germain is remembered once again, however briefly. Each time she is remembered, it comes with frustration and disbelief. Biographers in 1913 were bold in commenting on the importance of Germain's contributions -- ``\ldots all things considered, she was the most profoundly intellectual woman France has produced''~\cite{Mozans1913}. Jesse A. Fernandez Martinez wrote in 1946~\cite{Martinez1946}, 
\begin{quote}
``She showed great boldness in attacking a physical question which was at that time entirely outside the range of mathematical treatment and the more complicated cases of which had not been submitted to analysis at the time. \ldots It is a curious thing that a woman so deserving of recognition has not received it in a fuller measure.''
\end{quote}
Her work in number theory has been appreciated with far less controversy~\cite{Del2012}. Yet, while names like Euler, Lagrange, Laplace, Fourier, Navier, Cauchy, and Gauss are familiar to {\em all} students of mathematics, physics, and engineering, Germain's contributions remain obscure or misunderstood. Germain herself was modest and selfless about receiving credit for her work, 
\begin{quote}
``She rejoiced even when she saw her ideas made fruitful on occasion by other persons who adopted them. She repeatedly stated that it matters little who first arrives at an idea; rather what is significant is how far that idea can go. She said, happily, that her ideas had produced their fruit for science while not yielding anything for her reputation\ldots'' -- Guillaume Libri~\cite{Libri1832, Bucciarelli1980}
\end{quote}
We can and should do better for her today than her contemporaries.  Over three decades ago, Amy Dalm\`{e}dico was optimistic about how Germain's work would be perceived, noting that Germain ``\ldots believed, like the encyclopedists whom she had read, that her contributions to science would stand the tests of time and social prejudice on their own.''~\cite{Dalmedico1991}.

Although Euler is credited with identifying a means for calculating the principal curvatures of a surface, and thus what we would now recognize as the mean curvature, Germain recognized both the underlying mathematics of the mean curvature and its utility for addressing an important open problem in elasticity. I ask you to join me to referring to this measure as the ``Germain Curvature''.

\subsection{Acknowledgements}
I'm grateful for helpful feedback on early drafts of this essay from Howard A. Stone, Emma Lejeune, James Bird, Paul Barbone, and Meghan Gelardi Holmes.

\bibliography{GermainCurvature.bib}

\end{document}